\documentclass[11pt,a4paper]{article}

\usepackage{graphics,epic}
\usepackage{amsmath,amssymb, amsthm}
\usepackage[matrix, arrow]{xy}

\newtheorem{Thm}{Theorem}[section]
\newtheorem{Cor}[Thm]{Corollary}
\newtheorem{Lem}[Thm]{Lemma}

\theoremstyle{remark}

\newtheorem*{pf}{{\rm\textbf{Proof}}}

\begin{document}

\title{{\Large \bf On Defining Ideals or Subrings of Hall Algebras}\\
{\normalsize with an appendix by Andrew Hubery} }

\author{{\normalsize Dong Yang}\footnote{The author
acknowledges support by the AsiaLink network Algebras and
Representations in China and Europe, ASI/B7-301/98/679-11.}}

\date{}

\bigskip

\maketitle

\abstract{\begin{quote}Let $A$ be a finitary algebra over a finite
field $k$, and $A$-$mod$ the category of finite dimensional left
$A$-modules. Let $\mathcal{H}(A)$ be the corresponding Hall
algebra, and for a positive integer $r$ let $D_{r}(A)$ be the
subspace of $\mathcal{H}(A)$ which has a basis consisting of
isomorphism classes of modules in $A$-$mod$ with at least $r+1$
indecomposable direct summands. If $A$ is the path algebra of the
quiver of type $A_{n}$ with linear orientation, then $D_{r}(A)$ is
known to be the kernel of the map from the twisted Hall algebra to
the quantized Schur algebra indexed by $n+1$ and $r$. For any $A$,
we determine necessary and sufficient conditions for $D_{r}(A)$ to
be an ideal and some conditions for $D_{r}(A)$ to be a subring of
$\mathcal{H}(A)$. For $A$ the path algebra of a quiver, we also
determine
 necessary and sufficient conditions for $D_{r}(A)$ to be a subring of $\mathcal{H}(A)$.\end{quote}}

\section{Introduction}

Let $k$ be a finite field with $q$ elements and $A$ a $k$-algebra.
By an $A$-module we mean a finite dimensional left $A$-module.
Denote by $A$-$mod$ the category of $A$-modules. Assume $A$ is
finitary, i.e. $Ext^{1}(S_{1},S_{2})$ is a finite group for any
two (not necessarily different) simple objects in $A$-$mod$ ({\it
cf.}~\cite{R3}). Let $v=q^{1/2}$. Define the {\em Hall algebra}
$\mathcal{H}(A)$ to be the ${\mathbb{Z}}[v,v^{-1}]$-algebra with
basis the set of isomorphism classes $[X]$ of modules in $A$-$mod$
and with multiplication given by
\[\begin{array}{c} [M]\diamond [N] = \sum_{[X]}F^{X}_{M,N}[X] \end{array}\]
where $F^{X}_{M,N}$ is the number of submodules $U$ of $X$ such
that $U\cong N$ and $X/U\cong M$. For an $A$-module $M$, let
$s(M)$ be the number of indecomposable direct summands of $M$. For
an integer $r\geq 1$, let
\[
D_{r}(A)=\mathbb{Z}[v,v^{-1}]\{[M]\in \mathcal{H}(A) | s(M)\geq
r+1\}.\] For convenience we denote $D_{1}(A)$ by {\em $D(A)$}.

If in addition for any $A$-modules $M$, $N$, $Ext_{A}^{i}(M,N)=0$
for $i>>0$, one can define the {\em twisted Hall algebra
$\mathcal{H}_{*}(A)$} to be the ${\mathbb{Z}}[v,v^{-1}]$-algebra
with the same basis as $\mathcal{H}(A)$ and a twisted
multiplication
\[ [M] * [N] = v^{\langle M,N\rangle} [M]\diamond [N]\]
where $\langle M,N\rangle=\sum_{i\geq
0}(-1)^{i}dim_{k}Ext_{A}^{i}(M,N)$. One can check that $D_{r}(A)$
is an ideal (resp. subring) of $\mathcal{H}_{*}(A)$ if and only if
it is an ideal (resp. subring) of $\mathcal{H}(A)$. We refer
to~\cite{R2}~\cite{R3} for more on Hall algebras and twisted Hall
algebras.

For positive integers $n$ and $r$, let $S_{v}(n+1,r)$ be the
quantized Schur algebra of type $A_{n}$ and of degree $r$.
R.M.Green~\cite{RMG} shows that there is a map from the twisted
Hall algebra $\mathcal{H}_{*}(kL_{n})$ to $S_{v}(n+1,r)$ whose
kernel is exactly $D_{r}(kL_{n})$, where $L_{n}$ is the quiver of
type $A_{n}$ with linear orientation. In particular,
$D_{r}(kL_{n})$ is an ideal of $\mathcal{H}_{*}(kL_{n})$ (and
$\mathcal{H}(kL_{n})$) for all $r\geq 1$. This raises the
question, for which algebras $A$, $D_{r}(A)$ is an ideal of
$\mathcal{H}(A)$, or weaker, a subring? In general, we have

\begin{Thm} \label{T:MT1} The following conditions are equivalent,

{\rm(i)} $D(A)$ is an ideal of $\mathcal{H}(A)$,

{\rm(ii)} $A$ is serial, i.e. each indecomposable $A$-module is
 uniserial,

{\rm(iii)} $D_{r}(A)$ is an ideal of $\mathcal{H}(A)$ for all $r
\geq 1$,

{\rm(iv)} $D_{r}(A)$ is an ideal of $\mathcal{H}(A)$ for some
$r\geq 2$.
\end{Thm}

\begin{Thm} \label{T:MT2} Consider the following conditions,

{\rm(I)} Each indecomposable $A$-module has simple socle,

{\rm(I$^{\prime}$)} Each indecomposable $A$-module has simple top,
where the top of a module is defined as the quotient of the module
by its radical,

{\rm(II)} Each indecomposable $A$-module has simple top or simple
socle,

{\rm(III)} $D_{r}(A)$ is a subring of $\mathcal{H}(A)$ for all
$r\geq 1$,

{\rm(IV)} $D_{r}(A)$ is a subring of $\mathcal{H}(A)$ for some
$r\geq 2$,

{\rm(V)} $D(A)$ is a subring of $\mathcal{H}(A)$.\\
Then we have {\rm [(I) or (I)']$\iff$(III)$\iff$(IV)}
$\Longrightarrow$ {\rm (II)$\Longrightarrow$(V)}.
\end{Thm}

In the appendix by A.Hubery the equivalence {\rm (II)$\iff$(V)} is
proved for $A$ being a finite dimensional algebra. We conjecture
that this is true in general. \vskip5pt

Let $Q$ be a (finite) quiver. Let $s,t$ be the maps sending a path
to its starting vertex and terminating vertex respectively. Let
$A=kQ$ be the path algebra of $Q$ over $k$, where the product
$\alpha\beta$ of two paths $\alpha$ and $\beta$ of $Q$ is defined
as the composition of $\beta$ and $\alpha$ if
$t(\beta)=s(\alpha)$, and $0$ otherwise. Then $A$ is finitary and
$A$-$mod$ is equivalent to the category of finite dimensional
representations of $Q$. We refer to~\cite{ARS}~\cite{CB} for
representation theory of quivers. In the following, we will
identify an $A$-module with the corresponding representation of
$Q$. Let $\mathcal{H}(Q)=\mathcal{H}(A)$ be the Hall algebra,
$D(Q)=D(A)$ and $D_{r}(Q)=D_{r}(A)$. Then we have the following
theorem.

\begin{Thm} \label{T:quiver} {\rm(i)} $D_{r}(Q)$ is an ideal of $\mathcal{H}(Q)$ for all $r\geq 1$
if and only if $Q$ is a disjoint union of quivers of the form $L$
and $\Delta$, where

\hskip50pt {\setlength{\unitlength}{0.7pt}
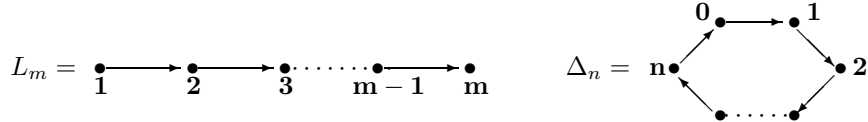
\begin{picture}(200,70)
\put(0,30){$\bullet$} \put(50,30){$\bullet$}
\put(100,30){$\bullet$} \put(150,30){$\bullet$}
\put(200,30){$\bullet$}

\put(7,34){$\vector(1,0){40}$} \put(57,34){$\vector(1,0){40}$}
\put(107,30){$\cdots\cdots\cdots$}
\put(157,34){$\vector(1,0){40}$}

\put(0,20){$\bf{1}$} \put(50,20){$\bf{2}$} \put(100,20){$\bf{3}$}
\put(140,20){$\bf{m-1}$} \put(200,20){$\bf{m}$}

\put(-50,30){ $L_{m}=$}

\end{picture}}
{\setlength{\unitlength}{0.7pt}
\begin{picture}(200,70)
\put(100,30){$\bullet$} \put(125,55){$\bullet$}
\put(165,55){$\bullet$} \put(190,30){$\bullet$}
\put(165,5){$\bullet$} \put(125,5){$\bullet$}

\put(107,37){$\vector(1,1){18}$} \put(132,59){$\vector(1,0){30}$}
\put(170,55){$\vector(1,-1){20}$}
\put(190,30){$\vector(-1,-1){20}$} \put(130,5){$\cdots\cdots$}
\put(124,11){$\vector(-1,1){18}$}

\put(115,60){$\bf{0}$} \put(175,60){$\bf{1}$}
\put(200,30){$\bf{2}$} \put(90,30){$\bf{n}$}

\put(40,30){ $\Delta_{n}=$}
\end{picture}}\\
$m$ is a positive integer and $n$ is a nonnegative integer (the
oriented cycle of $\Delta_{0}$ is a loop).

{\rm(ii)} $D(Q)$ is a subring of $\mathcal{H}(Q)$ if and only if
$Q$ is a disjoint union of quivers of the form  $L$, $\Delta$,
$V$, and $\Lambda$, where

{\setlength{\unitlength}{0.7pt}
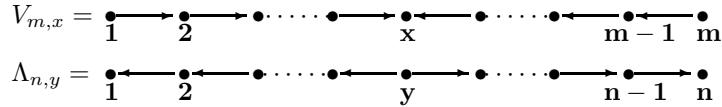
\begin{picture}(420,30)
\put(100,10){$\bullet$} \put(140,10){$\bullet$}
\put(180,10){$\bullet$} \put(220,10){$\bullet$}
\put(260,10){$\bullet$} \put(300,10){$\bullet$}
\put(340,10){$\bullet$} \put(380,10){$\bullet$}
\put(420,10){$\bullet$}

\put(107,14){$\vector(1,0){30}$} \put(147,14){$\vector(1,0){30}$}
\put(187,10){$\cdots\cdots$} \put(227,14){$\vector(1,0){30}$}
\put(298,14){$\vector(-1,0){30}$} \put(305,10){$\cdots\cdots$}
\put(378,14){$\vector(-1,0){30}$}
\put(418,14){$\vector(-1,0){30}$}

\put(100,0){$\bf{1}$} \put(140,0){$\bf{2}$} 
\put(260,0){$\bf{x}$} 
\put(370,0){$\bf{m-1}$} \put(420,0){$\bf{m}$}

\put(50,10){$V_{m,x}=$}
\end{picture}}

{\setlength{\unitlength}{0.7pt}
\begin{picture}(400,30)
\put(100,10){$\bullet$} \put(140,10){$\bullet$}
\put(180,10){$\bullet$} \put(220,10){$\bullet$}
\put(260,10){$\bullet$} \put(300,10){$\bullet$}
\put(340,10){$\bullet$} \put(380,10){$\bullet$}
\put(420,10){$\bullet$}

\put(138,14){$\vector(-1,0){30}$}
\put(178,14){$\vector(-1,0){30}$} \put(185,10){$\cdots\cdots$}
\put(258,14){$\vector(-1,0){30}$} \put(267,14){$\vector(1,0){30}$}
\put(305,10){$\cdots\cdots$} \put(347,14){$\vector(1,0){30}$}
\put(387,14){$\vector(1,0){30}$}

\put(100,0){$\bf{1}$} \put(140,0){$\bf{2}$} \put(260,0){$\bf{y}$}
\put(370,0){$\bf{n-1}$} \put(420,0){$\bf{n}$}

\put(50,10){$\Lambda_{n,y}=$}
\end{picture}}\\
$m$, $n\geq 3$ are positive integers, $x\in\{2,\cdots,m-1\}$,
$y\in\{2,\cdots,n-1\}$.

{\rm(iii)} The following conditions are equivalent,

\ \ \ \ {\rm(a)} $D_{r}(Q)$ is a subring of $\mathcal{H}(Q)$ for
all $r\geq 1$,

\ \ \ \ {\rm(b)} $D_{r}(Q)$ is a subring of $\mathcal{H}(Q)$ for
some $r\geq 2$,

\ \ \ \ {\rm(c)} $Q$ is a disjoint union of quivers of the form
$L$, $\Delta$, and $V$, or a disjoint union of quivers of the form
$L$, $\Delta$, and $\Lambda$.
\end{Thm}

When $\mathcal{H}(Q)$ is replaced by $\mathcal{H}^{nil}(Q)$, the
subalgebra of $\mathcal{H}(Q)$ with basis the isomorphism classes
of finite dimensional nilpotent representations of $Q$, and
$D_{r}(Q)$ by $D^{nil}_{r}(Q)=D_{r}(Q)\cap \mathcal{H}^{nil}(Q)$,
Theorem~\ref{T:quiver} still holds. \vskip5pt

\noindent{\bf Acknowledgement} The author thanks Andrew Hubery for
writing the appendix. He is also grateful to the referee for very
helpful remarks which make the proofs of the main theorems more
transparent.

\section{Proof of the theorems}

Assume $0\rightarrow N\rightarrow X\rightarrow M\rightarrow 0$ is
an exact sequence in $A$-$mod$, i.e. $[X]$ is a summand (up to
scalar) of $[M]\diamond [N]$, then the socle $soc(N)$ of $N$ and
the top $top(M)$ of $M$ are direct summands of $soc(X)$ and
$top(X)$ respectively.

\begin{Lem}\label{L:left} The following conditions are equivalent,

{\rm(I)} Each indecomposable $A$-module has simple socle,

{\rm(VI)} $D_{r}(A)$ is a left ideal of $\mathcal{H}(A)$ for all
$r\geq 1$,

{\rm (VII)} $D_{r}(A)$ is a left ideal of $\mathcal{H}(A)$ for
some $r\geq 1$.

\end{Lem}
\begin{pf} (I)$\Rightarrow$(VI) Assume $r\geq 1$, and $[M]\in D_{r}(A)$. Then by (I)
 $soc(M)$ has at least $r+1$ direct summands.
Therefore if $X$ is an extension of some $A$-module $N$ by $M$,
then $soc(X)$ has at least $r+1$ direct summands, and hence $X$
has at least $r+1$ direct summands by (I). Therefore $D_{r}(A)$ is
a left ideal of $\mathcal{H}(A)$.

(VII)$\Rightarrow$(I) Suppose on the contrary that there exists an
indecomposable module $M$ such that $soc(M)$ is decomposable. Let
$r\geq 1$ be any integer and $N$ an indecomposable $A$-module.
Then $0\rightarrow soc(M)\oplus N^{\oplus r-1}\rightarrow M\oplus
N^{\oplus r-1}\rightarrow M/soc(M)\rightarrow 0$ is an exact
sequence, and $soc(M)\oplus N^{\oplus r-1}$ has at least $r+1$
direct summands but $M\oplus N^{\oplus r-1}$ has exactly $r$
indecomposable direct summands. Thus $D_{r}(A)$ is not a left
ideal of $\mathcal{H}(A)$ for all $r\geq 1$, contradicting (VII).
$\square$

\end{pf}

Dually, we have

\begin{Lem}\label{L:right} The following conditions are equivalent,

{\rm(I$^{\prime}$)} Each indecomposable $A$-module has simple top,

{\rm(VI$^{\prime}$)} $D_{r}(A)$ is a right ideal of
$\mathcal{H}(A)$ for all $r\geq 1$,

{\rm (VII$^{\prime}$)} $D_{r}(A)$ is a right ideal of
$\mathcal{H}(A)$ for some $r\geq 1$.
\end{Lem}
\vskip5pt

\noindent{\bf Proof of Theorem~\ref{T:MT1} : }

It follows from Lemma~\ref{L:left} and Lemma~\ref{L:right} since a
module which has simple top and simple socle is uniserial.
$\square$

\vskip10pt \noindent{\bf Proof of Theorem~\ref{T:MT2} : }

[(I) or (I$^{\prime}$)] $\Rightarrow$(III) It follows from
Lemma~\ref{L:left} and Lemma~\ref{L:right}.

(IV)$\Rightarrow$ [(I) or (I$^{\prime}$)] Suppose $M$, $N$ are two
indecomposable $A$-modules such that $soc(M)$ is not simple and
$top(N)$ is not simple. Let $r\geq 2$ be any integer. Then
$0\rightarrow soc(M)\oplus rad(N)^{\oplus r-1}\rightarrow M\oplus
N^{\oplus r-1} \rightarrow M/soc(M)\oplus top(N)^{\oplus
r-1}\rightarrow 0$ is an exact sequence, and $M/soc(M)\oplus
top(N)^{\oplus r-1}$ and $soc(M)\oplus rad(N)^{\oplus r-1}$ both
have at least $r+1$ direct summands, but $M\oplus N^{\oplus r-1}$
has exactly $r$ indecomposable direct summands. Thus $D_{r}(A)$ is
not a subring of $\mathcal{H}(A)$, contradicting (IV).

(II)$\Rightarrow$(V) Let $M$, $N$ be two decomposable $A$-modules.
In particular, $top(M)$ and $soc(N)$ are decomposable modules. Let
$X$ be an extension of $M$ by $N$. Then both $top(X)$ and $soc(X)$
are decomposable. Since by (II) each indecomposable $A$-module has
simple top or simple socle, $X$ is decomposable. $\square$

\begin{Cor} Let $r\geq 2$ be an integer, then $D_{r}(A)$ is a subring
of $\mathcal{H}(A)$ if and only if $D_{r}(A)$ is a left ideal or a
right ideal of $\mathcal{H}(A)$.
\end{Cor}

\vskip5pt

In the sequel we fix our attention on quivers. We observe the
following facts.

Assume $Q$, $Q'$ are two quivers and there is an exact functor
$F:kQ'$-$mod\rightarrow kQ$-$mod$ preserving indecomposability of
representations. Then for any positive integer $r$ the space
$D_{r}(Q')$ is not a subring (resp. an ideal) of $\mathcal{H}(Q')$
implies that $D_{r}(Q)$ is not a subring (resp. an ideal) of
$\mathcal{H}(Q)$. The statement is true if we replace $D_{r}$ by
$D_{r}^{nil}$ and $\mathcal{H}$ by $\mathcal{H}^{nil}$ provided
$F$ sends nilpotent representations to nilpotent representations.
In particular,

(i) Let $Q'$ be a subquiver of $Q$, i.e. the sets of vertices and
arrows of $Q'$ are subsets of those of $Q$. Then a representation
of $Q'$ can be regarded as a representation of $Q$ by putting zero
vector spaces and zero maps on vertices and arrows of $Q$ which
are different from those of $Q'$ respectively. This defines a
(covariant) exact functor $F:kQ'$-$mod\rightarrow kQ$-$mod$ which
preserves the indecomposability of representations and sends
nilpotent representations to nilpotent representations.

(ii) Let $Q$ be a quiver and $Q^{op}$ the quiver with the same
underlying diagram as $Q$ but opposite orientation. To a
representation of $Q$ we associate a representation of $Q^{op}$ by
taking the dual of all vectors spaces and all maps. This defines a
(contravariant) exact functor from $kQ$-$mod$ to $kQ^{op}$-$mod$
which preserves indecomposability and sends nilpotent
representations to nilpotent representations. In particular,
$D_{r}(Q)$ is a subring (resp. an ideal) of $\mathcal{H}(Q)$ if
and only if $D_{r}(Q^{op})$ is a subring (resp. an ideal) of
$\mathcal{H}(Q^{op})$. If we replace $\mathcal{H}$ by
$\mathcal{H}^{nil}$ and $D_{r}$ by $D_{r}^{nil}$, the statement is
true.

\begin{Thm} \label{T:connectedquiver} Let $Q$ be a connected quiver, and $L$, $\Delta$, $V$ and $\Lambda$ as in Theorem~\ref{T:quiver}.

{\rm(i)} The following conditions are equivalent,

\ \ \ \ {\rm(a)} $D(Q)$ is an ideal of $\mathcal{H}(Q)$,

\ \ \ \ {\rm(b)} $D_{r}(Q)$ is an ideal of $\mathcal{H}(Q)$ for
all $r\geq 1$,

\ \ \ \ {\rm(c)} $Q$ is of the form $L$ or $\Delta$.

{\rm(ii)} The following conditions are equivalent,

\ \ \ \ {\rm(a)} $D(Q)$ is a subring of $\mathcal{H}(Q)$,

\ \ \ \ {\rm(b)} $D_{r}(Q)$ is a subring of $\mathcal{H}(Q)$ for
all $r\geq 1$,

\ \ \ \ {\rm(c)} $Q$ is of the form $L$, $\Delta$, $V$, or
$\Lambda$.
\end{Thm}

To prove Theorem~\ref{T:connectedquiver}, we need the following
Lemmas~\ref{L:An}-~\ref{L:twoloops}.\vskip10pt

\begin{Lem} \label{L:An} Let $Q_{1}=$ {\setlength{\unitlength}{0.7pt}
\begin{picture}(180,40)
\put(10,0){$\bullet$} \put(40,30){$\bullet$} \put(70,0){$\bullet$}

\put(100,-30){$\bullet$} \put(160,-30){$\bullet$}
\put(130,-60){$\bullet$}

\put(38,29){$\vector(-1,-1){22}$} \put(47,29){$\vector(1,-1){22}$}
\put(70,-10){$\ddots$} \put(84,-23){$\ddots$}
\put(108,-30){$\vector(1,-1){22}$}
\put(158,-30){$\vector(-1,-1){22}$}

\put(0,0){$\bf{1}$} \put(40,40){$\bf{2}$}
\put(130,-70){$\bf{n-1}$} \put(170,-30){$\bf{n}$}

\end{picture}}.\\
$\bf{}$\\
$\bf{}$\\
$\bf{}$\\
Then $D(Q_{1})$ is not a subring of $\mathcal{H}(Q_{1})$.
\end{Lem}

\begin{pf} Let, for a moment, $Q$ be any quiver. To prove
$D(Q)$ is not a subring of $\mathcal{H}(Q)$, it suffices to prove
that there exists an indecomposable $kQ$-module which has a
decomposable submodule with quotient also decomposable. For
$Q_{1}$, the indecomposable module with dimension vector
$(1,1,\cdots,1)$ has a submodule with dimension vector
$(1,0,1,\cdots,1,0)$ with quotient having dimension vector
$(0,1,0,\cdots,0,1)$. $\square$
\end{pf}

\begin{Lem} \label{L:D4} Let $Q_{2}$ be a quiver of type $D_{4}$, i.e. the underlying graph of $Q_{2}$ is  {\setlength{\unitlength}{0.7pt}
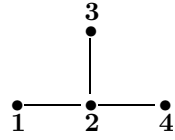
\begin{picture}(120,40)
\put(0,-20){$\bullet$} \put(40,-20){$\bullet$}
\put(80,-20){$\bullet$} \put(40,20){$\bullet$}

\put(7.5,-16){$\line(1,0){30}$} \put(47.5,-16){$\line(1,0){30}$}
\put(43,-10){$\line(0,1){30}$}

\put(0,-30){$\bf{1}$} \put(40,-30){$\bf{2}$}
\put(80,-30){$\bf{4}$} \put(40,30){$\bf{3}$}
\end{picture} }\\$\bf{}$\\
Then $D(Q_{2})$ is not a subring of $\mathcal{H}(Q_{2})$.
\end{Lem}

\begin{pf} Let $M$ be the
indecomposable module with dimension vector $(1,2,1,1)$. Then
either both $rad(M)$ and $M/rad(M)$ are decomposable, or both
$soc(M)$ and $M/soc(M)$ are decomposable. $\square$
\end{pf}

\begin{Lem} \label{L:kronecker} Let $Q_{3}=$ {\setlength{\unitlength}{0.7pt}
\begin{picture}(80,15)
\put(10,0){$\bullet$} \put(60,0){$\bullet$}

\put(20,6.5){$\vector(1,0){37}$} \put(20,1.5){$\vector(1,0){37}$}

\put(0,0){$\bf{1}$} \put(70,0){$\bf{2}$}
\end{picture}} be the Kronecker quiver.
Then $D(Q_{3})$ is not a subring of $\mathcal{H}(Q_{3})$.
\end{Lem}
\begin{pf} Consider the indecomposable module $M$ with
dimension vector $(3,2)$. Then both $soc(M)$ and $M/soc(M)$ are
decomposable. $\square$
\end{pf}

\begin{Lem} \label{L:shoot} Let \ $Q_{4}=$ {\setlength{\unitlength}{0.7pt}
\begin{picture}(110,20)
\put(0,5){$\bullet$} \put(50,5){$\bullet$} \put(100,5){$\bullet$}

\put(10,11){$\vector(1,0){35}$} \put(45,6){$\vector(-1,0){35}$}
\put(60,9){$\vector(1,0){35}$}

\put(0,-5){$\bf{1}$} \put(50,-5){$\bf{2}$} \put(100,-5){$\bf{3}$}

\end{picture}} , \ $Q_{5}=$ {\setlength{\unitlength}{0.7pt}
\begin{picture}(110,20)
\put(0,5){$\bullet$} \put(50,5){$\bullet$} \put(100,5){$\bullet$}

\put(10,11){$\vector(1,0){35}$} \put(45,6){$\vector(-1,0){35}$}
\put(95,9){$\vector(-1,0){35}$}

\put(0,-5){$\bf{1}$} \put(50,-5){$\bf{2}$} \put(100,-5){$\bf{3}$}

\end{picture}} , \ $Q_{6}$= {\setlength{\unitlength}{0.7pt}
\begin{picture}(90,30)
\qbezier(20,5)(10,-10)(2,10) \qbezier(20,15)(10,30)(2,10)

\put(23,6){$\bullet$} \put(68,10){$\vector(-1,0){35}$}
\put(73,6){$\bullet$}

\drawline(20,15)(20,20) \drawline(20,15)(15,17)

\put(23, -5){$\bf{1}$} \put(73,-5){$\bf 2$}

\end{picture}}, \ $Q_{7}$= {\setlength{\unitlength}{0.7pt}
\begin{picture}(90,30)
\qbezier(20,5)(10,-10)(2,10) \qbezier(20,15)(10,30)(2,10)

\put(23,6){$\bullet$} \put(33,10){$\vector(1,0){35}$}
\put(73,6){$\bullet$}

\drawline(20,15)(20,20) \drawline(20,15)(15,17)

\put(23,-5){$\bf 1$} \put(73,-5){$\bf 2$}

\end{picture}}.
Then $D(Q_{i})$ is not a subring of $\mathcal{H}(Q_{i})$,
$i=4,5,6,7$.
\end{Lem}
\begin{pf} Sending {\setlength{\unitlength}{0.7pt}
\begin{picture}(100,25)
\qbezier(30,-1)(20,-16)(12,4) \qbezier(30,9)(20,24)(12,4)
\drawline(30,9)(30,14) \drawline(30,9)(25,11)
 \put(33,0){$V_{1}$}
\put(50,-1){$\longleftarrow$} \put(78,0){$V_{2}$}
\put(0,9){\footnotesize $f_{1}$} \put(53,9){\footnotesize $f_{2}$}
\end{picture}}  to   {\setlength{\unitlength}{0.7pt}
\begin{picture}(120,25)
\put(0,0){ $V_{1}$ } \put(25,-3){$\longleftarrow$}
\put(25,3){$\longrightarrow$}

\put(50,0){ $V_{1}$ $\longleftarrow$ $V_{2}$}

\put(30,10){\footnotesize $f_{1}$} \put(30,-10){\footnotesize
$id$} \put(80,10){\footnotesize $f_{2}$}

\end{picture}}  yields an exact functor from $kQ_{6}$-$mod$
to $kQ_{5}$-$mod$ which preserves the indecomposability of
representations. Moreover, $Q_{4}$ and $Q_{5}$ are opposite to
each other and so are $Q_{6}$ and $Q_{7}$. Therefore it is enough
to prove the statement for $Q_{6}$. We have
\[
{\setlength{\unitlength}{0.7pt}
\begin{picture}(400,25)
\qbezier(40,10)(30,-5)(22,15) \qbezier(40,20)(30,35)(22,15)
\drawline(40,20)(40,25) \drawline(40,20)(35,22)
 \put(43,11){$k$}
\put(55,10){$\longleftarrow$} \put(82,11){$k$} \put(0,10){ [ }
\put(17,20){\footnotesize $0$} \put(63,15){\footnotesize $0$}
\put(92,10){ ] $\diamond$ [ }

\qbezier(151,10)(141,-5)(133,15) \qbezier(151,20)(141,35)(133,15)
\drawline(151,20)(151,25) \drawline(151,20)(146,22)
 \put(154,11){$k^{3}$}
\put(169,10){$\longleftarrow$} \put(196,11){$0$} \put(128,20){$f$}
\put(206,10){ ] $=$  [ }

\qbezier(272,10)(262,-5)(254,15) \qbezier(272,20)(262,35)(254,15)
\drawline(272,20)(272,25) \drawline(272,20)(267,22)
 \put(275,11){$k^{4}$}
\put(290,10){$\longleftarrow$} \put(317,11){$k$}
\put(245,21){$g_{1}$} \put(298,20){$g_{2}$} \put(327,10){ ] +
others }
\end{picture}}\]
where  $f = \footnotesize \left ( \begin{array}{ccc}
                                     0 & 0 & 0\\
                                     1 & 0 & 0\\
                                     0 & 0 & 0\end{array}\right
                                     )$,
       $g_{1} = \footnotesize \left ( \begin{array}{cccc}
                                     0 & 0 & 0 & 0\\
                                     1 & 0 & 0 & 0\\
                                     0 & 1 & 0 & 0\\
                                     0 & 0 & 0 & 0\end{array}\right )$,
      $g_{2} = \footnotesize \left ( \begin{array}{c}
                                      0\\
                                      1\\
                                      0\\
                                      1\end{array}\right )$.

Both modules on the left hand side are decomposable. Let $M$
denote the module on the right hand side. Then
\[
End(M) = \Big\{ (A,B) \Big| A= \footnotesize \left(
\begin{array}{cccc}
                                            a & 0 & 0 & 0 \\
                                            b & a & 0 & 0 \\
                                            c & b & a & -b\\
                                            d & 0 & 0 & a\end{array}\right ),
                      B=a, {\rm \ where }\ a,b,c,d\in k \Big\}.\] It has a unique maximal
ideal $\{(A,B)| a=0 \}$, so it is local. Therefore, $M$ is
indecomposable. $\square$

\end{pf}

\begin{Lem} \label{L:twoloops} Let $Q_{8}$=
{\setlength{\unitlength}{0.7pt}
\begin{picture}(60,20)
\qbezier(20,0)(10,-15)(2,5) \qbezier(20,10)(10,25)(2,5)

\put(20,1){$\bullet$}

\drawline(20,10)(20,15) \drawline(20,10)(15,12)

\qbezier(26,0)(36,-15)(44,5) \qbezier(26,10)(36,25)(44,5)

\drawline(26,10)(26,15) \drawline(26,10)(31,12)

\put(20,-13){$\bf 1$}
\end{picture}}.
Then $D(Q_{8})$ is not a subring of $\mathcal{H}(Q_{8})$.
\end{Lem}
\begin{pf} We have
\[
{\setlength{\unitlength}{0.7pt}
\begin{picture}(350,30)
\put(5,15){ [ } \qbezier(40,15)(30,0)(22,20)
\qbezier(40,25)(30,40)(22,20)

\drawline(40,25)(40,30) \drawline(40,25)(35,27)

\put(43,16){$k^{2}$}

\qbezier(57,15)(67,0)(75,20) \qbezier(57,25)(67,40)(75,20)

\drawline(57,25)(57,30) \drawline(57,25)(62,27)

\put(17,24){\footnotesize $0$} \put(78,24){\footnotesize $0$}

\put(90,15){ ] $\diamond$ [ }

\qbezier(150,15)(140,0)(132,20) \qbezier(150,25)(140,40)(132,20)

\drawline(150,25)(150,30) \drawline(150,25)(145,27)

\put(153,16){$k^{2}$}

\qbezier(167,15)(177,0)(185,20) \qbezier(167,25)(177,40)(185,20)

\drawline(167,25)(167,30) \drawline(167,25)(172,27)

\put(127,24){\footnotesize $0$} \put(188,24){\footnotesize $0$}

\put(200,15){ ] $=$ [ }

\qbezier(265,15)(255,0)(247,20) \qbezier(265,25)(255,40)(247,20)

\drawline(265,25)(265,30) \drawline(265,25)(260,27)

\put(268,16){$k^{4}$}

\qbezier(282,15)(292,0)(300,20) \qbezier(282,25)(292,40)(300,20)

\drawline(282,25)(282,30) \drawline(282,25)(287,27)

\put(238,24){$g_{1}$} \put(303,24){$g_{2}$}

\put(318,15){ ] + others }

\end{picture}}\]
where $g_{1} = \footnotesize \left ( \begin{array}{cccc}
                                     0 & 0 & 0 & 0\\
                                     0 & 0 & 0 & 0\\
                                     1 & 0 & 0 & 0\\
                                     0 & 1 & 0 & 0\end{array}\right )$,
      $g_{2} = \footnotesize \left ( \begin{array}{cccc}
                                      0 & 0 & 0 & 0\\
                                      0 & 0 & 0 & 0\\
                                      0 & 0 & 0 & 0\\
                                      1 & 0 & 0 & 0\end{array}\right )$.

Both modules on the left hand side are decomposable. Let $M$
denote the module on the right hand side. Then
\[End(M) = \Big\{ A
\Big| A=\footnotesize \left(
\begin{array}{cccc}
                                            a & 0 & 0 & 0 \\
                                            b & a & 0 & 0 \\
                                            d & c & a & 0\\
                                            e & f & b & a\end{array}\right ),
                  {\rm \ where\ } a,b,c,d,e,f\in k\Big\}.\]
It has a unique maximal ideal $\{ A | a=0 \}$, so it is local.
Hence $M$ is indecomposable. $\square$
\end{pf}

\vskip10pt \noindent {\bf Proof of Theorem~\ref{T:connectedquiver}
: }

Let $Q$ be a connected quiver. It is known that $kQ$ is serial if
and only if $Q$ is of the form $L$ or $\Delta$. Thus (i) follows
from Theorem~\ref{T:MT1}. Now let us prove (ii).

If $Q$ is of the form $L$, $\Delta$, $V$ or $\Lambda$, then
Theorem~\ref{T:MT2} [(I) or (I$^{\prime}$)] holds, and hence
$D_{r}(Q)$ is a subring of $\mathcal{H}(Q)$ for all $r\geq 1$. By
Theorem~\ref{T:MT2} it remains to show that if $Q$ is not of the
form $L$, $\Delta$, $V$ or $\Lambda$ then $D(Q)$ is not a subring
of $\mathcal{H}(Q)$. We prove case by case.

Case 1. $Q$ is of type $A_{n}$ but not of the form $L$, $V$, or
$\Lambda$. Then $Q$ has a subquiver of the form $Q_{1}$. Therefore
by Lemma~\ref{L:An} we have that $D(Q)$ is not a subring of
$\mathcal{H}(Q)$.

Case 2. $Q$ is of type $\tilde{A}_{n}$ but not of the form
$\Delta$. Then there exists an exact functor from $kQ_{3}$-$mod$
to $kQ$-$mod$ which preserves indecomposability of representations
({\it cf.}~\cite{DR}~\cite{R1}). Therefore by
Lemma~\ref{L:kronecker} $D(Q)$ is not a subring of
$\mathcal{H}(Q)$.

Case 3. $Q$ has a proper subquiver which is of the form $\Delta$.
Then $Q$ has a double-loop or a subquiver of type $D_{4}$, or a
subquiver of the form $Q_{4}$, $Q_{5}$, $Q_{6}$, or $Q_{7}$.
Therefore it follows from Lemma ~\ref{L:twoloops},
Lemma~\ref{L:D4} and Lemma~\ref{L:shoot} that $D(Q)$ is not a
subring of $\mathcal{H}(Q)$.

Case 4. Otherwise. Then $Q$ has a subquiver of type $D_{4}$ or a
Kronecker subquiver. Thus by Lemma~\ref{L:D4} and
Lemma~\ref{L:kronecker} $D(Q)$ is not a subring of
$\mathcal{H}(Q)$. $\square$

\vskip10pt \noindent{\bf Proof of Theorem~\ref{T:quiver} : }

(i) It follows, by Theorem~\ref{T:MT1}, from the fact that $kQ$ is
serial if and only if $Q$ is a disjoint union of quivers of the
form $L$ and $\Delta$.

(ii) Assume $D(Q)$ is a subring of $\mathcal{H}(Q)$. It follows
from Theorem~\ref{T:connectedquiver} that $Q$ is a disjoint union
of quivers of the form $L$, $\Delta$, $V$ and $\Lambda$. If $Q$ is
such a quiver, then each $kQ$-module has simple top or simple
socle. By Theorem~\ref{T:MT2}, $D(Q)$ is a subring of
$\mathcal{H}(Q)$.

(iii) Assume $Q$ is a disjoint union of quivers of the form $L$,
$\Delta$, $V$ and $\Lambda$. Then each $kQ$-module has simple
socle if and only if $Q$ is a disjoint union of quivers of the
form $L$, $\Delta$ and $V$; each $kQ$-module has simple top if and
only if $Q$ is a disjoint union of quivers of the form $L$,
$\Delta$, and $\Lambda$. The desired result follows
fromTheorem~\ref{T:MT2}. $\square$ \vskip5pt

We can follow the same procedure to prove

\begin{Thm} {\rm(i)} $D_{r}^{nil}(Q)$ is an ideal of $\mathcal{H}^{nil}(Q)$ for any integer $r\geq 1$
if and only if $Q$ is a disjoint union of quivers of the form $L$
and $\Delta$.

{\rm(ii)} $D^{nil}(Q)$ is a subring of $\mathcal{H}^{nil}(Q)$ if
and only if $Q$ is a disjoint union of quivers of the form  $L$,
$\Delta$, $V$, and $\Lambda$.

 {\rm(iii)} The following conditions are
equivalent,

\ \ \ \ {\rm(a)} $D_{r}^{nil}(Q)$ is a subring of
$\mathcal{H}^{nil}(Q)$ for any integer $r\geq 1$,

\ \ \ \ {\rm(b)} $D_{r}^{nil}(Q)$ is a subring of
$\mathcal{H}^{nil}(Q)$ for some integer $r\geq 2$,

\ \ \ \ {\rm(c)} $Q$ is a disjoint union of quivers of the form
$L$, $\Delta$, and $V$, or a disjoint union of quivers of the form
$L$, $\Delta$, and $\Lambda$.
\end{Thm}

\vskip20pt

\section{Appendix by Andrew Hubery}

Let $k$ be a finite field and $A$ a finite dimensional
$k$-algebra. Let $\mathcal H(A)$ denote the Ringel-Hall algebra of
$A$ and define $D(A)$ to the be subspace of all decomposable
modules.

\begin{Thm}\label{T:Hub}
$D(A)$ is a subring of $\mathcal H(A)$ if and only if every
indecomposable $A$-module has either simple top or simple socle.
\end{Thm}

We shall use the following characterisation, due to Tachikawa
\cite{Tach}.
\begin{Thm}
Every indecomposable $A$-module has simple top or simple socle if
and only if
\begin{enumerate}
\item every indecomposable projective module has radical a sum of
at most two uniserial modules (and dually for indecomposable
injective modules), and \item if an indecomposable projective has
decomposable socle, then the injective envelopes of these simples
are uniserial (and dually for injectives).
\end{enumerate}
\end{Thm}

\begin{proof}
It is clear that if every indecomposable has simple top or simple
socle, then $D(A)$ is a subring, since every extension of
decomposable modules must remain decomposable. Let us therefore
assume that $D(A)$ is a subring, and write $Q$ for the valued
quiver of $A$.

We first consider those indecomposable modules of Loewy length
two. There are no valued arrows with valuation $(a,b)$ for
$ab\geq3$, so in particular, there is no Kronecker subquiver, and
hence no vertex with a double loop. Also, there are at most two
arrows starting at each vertex, and if there are two such arrows,
then they are both unvalued. Dually for arrows ending at a given
vertex. Finally, we can have no subquiver of the form
$\cdot\leftarrow\cdot\rightarrow\cdot\leftarrow$.

We now consider indecomposable modules of Loewy length three.
Suppose we have a vertex $i$ with an arrow $\alpha$ ending at $i$
and two arrows $\beta$ and $\gamma$ starting at $i$. We know from
the above that both $\beta$ and $\gamma$ are unvalued, but
$\alpha$ may be valued. We show that there can be no sincere
indecomposable module for this subquiver.

It is sufficient to consider the valued graph
$$\xymatrix@R=0pt{&&2\\1\ar[r]^{(a,b)} &i\ar[ur]\ar[dr]&&ab\leq2.\\&&3}$$
The corresponding $k$-species $\Lambda$ is given by
$$\Lambda=\begin{pmatrix}F&0&0&0\\H&G&0&0\\H&H&G&0\\H&H&0&G\end{pmatrix},$$
where $F/k$, $G/k$ and $H/k$ are field extensions of degrees $a$,
$b$ and $ab$ respectively, and ${}_GH_F$ has the natural bimodule
structure. Note that as $ab\leq2$, then either $a=1$, so $F=k$ and
$G=H$, or $b=1$, so $F=H$ and $G=k$. Now,
$\mathrm{rad}^2\Lambda\cong H\oplus H$ and hence the only possible
relations are the zero relations $\beta\alpha=0$ or
$\gamma\alpha=0$. Thus, if there exists a sincere indecomposable
module, then there are no relations and $\mathrm{mod}\,\Lambda$
embeds into $\mathrm{mod}\,A$. This gives a contradiction, since
$D(\Lambda)$ is not a subring of $\mathcal H(\Lambda)$.

The dual argument works whenever there are two arrows ending at
$i$ and an arrow starting at $i$.

It follows that for each indecomposable projective $P$,
$\mathrm{rad}\,P$ is the sum of at most two modules. These modules
must be uniserial, since if not, then we are in the situation
above for some vertex $i$: that is, there exists an arrow ending
at $i$, two arrows starting at $i$ and no zero relations, a
contradiction.

Now let $P$ be an indecomposable projective module such that
$\mathrm{soc}\,P$ is decomposable. Write $\mathrm{rad}\,P=U_1+U_2$
as a sum of two uniserial modules and let $j$ be the vertex
corresponding to $S=\mathrm{soc}\,U_1$. Suppose that $I(S)$ is not
uniserial. Then the module $M=P/\mathrm{rad}\,U_2$ is
indecomposable and we claim that there exists an (unvalued) arrow
$\alpha:i\to j$ such that $\alpha\cdot M=0$.

If the socles of $U_1$ and $U_2$ are non-isomorphic, this is
clear, so suppose that $\mathrm{soc}\,P\cong S^2$. Then $P$ has
Loewy length at least three and we may assume that $U_2$ has Loewy
length at least two. This proves the claim.

Now, there is a natural non-split extension of $S_i$ by $M$
yielding an indecomposable with decomposable radical and
decomposable top, a contradiction.

We clearly have the dual statements involving indecomposable
injective modules, and hence the conditions of Tachikawa's Theorem
are fulfilled.
\end{proof}

As a corollary, we extend Theorem\ref{T:connectedquiver} to all
hereditary algebras.
\begin{Cor}
Let $A$ be a connected hereditary $k$-algebra. Then $D(A)$ is a
subring if and only if the quiver of $A$ is either an oriented
cycle, of type $\mathbb A$ and having either a unique sink or a
unique source, or of type $\mathbb B$ or $\mathbb C$ with a linear
orientation.
\end{Cor}

\begin{proof}
Assume that $A$ is not an oriented cycle. We know from \cite{DR}
that if $A$ is of tame representation type, then there exists an
embedding into $\mathrm{mod}\,A$ of the module category for some
tame bimodule. This has valuation $(a,b)$ with $ab=4$, a
contradiction. Thus $A$ must be representation finite. Now, by the
previous arguments and Tachikawa's criteria, $A$ must be of type
$\mathbb A$, $\mathbb B$ or $\mathbb C$ with the required
restrictions on the the orientation.
\end{proof}

\bibliographystyle{amsplain}
\bibliography{xbib}

Dong Yang

Department of Mathematical Sciences, Tsinghua University,
Beijing100084, P.R.China.

Department of Mathematics, University of Leicester, University
Road, Leicester LE1 7RH, United Kingdom.

{\it Email address : }yangdong98@mails.tsinghua.edu.cn . \vskip5pt

Andrew Hubery

Universit\"{a}t Paderborn, Institut f\"{u}r Mathematik, Warburger
Stra{\ss}e, 100 33098 Paderborn, Germany.

{\it Email address : }hubery@math.uni-paderborn.de .

\end{document}